# SUPER-MATHEMATICS FUNCTIONS

by Mircea Eugen Şelariu
Polytechnic University of Timişoara, Romania

[translated from Romanian by Marian Niţu and Florentin Smarandache]

In this paper we talk about the so-called **Super-Mathematics Functions (SMF)**, which often constitute the base for generating technical, neo-geometrical, therefore less artistic objects.

These functions are the results of 38 years of research, which began at University of Stuttgart in 1969. Since then, 42 related works have been published, written by over 19 authors, as shown in the References.

The name was given by the regretted mathematician Professor Emeritus Doctor Engineer **Gheorghe Silas** who, at the presentation of the very first work in this domain, during the First National Conference of Vibrations in Machine Constructions, Timişoara, Romania, 1978, named CIRCULAR EX-CENTRIC FUNCTIONS, declared: "Young man, you just discovered not only "some functions, but a new mathematics, a **supermathematics**!" I was glad, at my age of 40, like a teenager. And I proudly found that he might be right!

The prefix **super** is justified today, to point out the birth of the new complements in mathematics, joined together under the name of Ex-centric Mathematics (EM), with much more important and infinitely more numerous entities than the existing ones in the **actual mathematics,** which we are obliged to call it **Centric Mathematics (CM.)**

To each entity from CM corresponds an infinity of similar entities in EM, therefore the Supermathematics (SM) is the reunion of the two domains: **SM = CM ∪ EM,** where **CM** is a particular case of null ex-centricity of **EM**. Namely, **CM = SM(e = 0)**. To each known function in **CM** corresponds an infinite family of functions in **EM**, and in addition, a series of new functions appear, with a wide range of applications in mathematics and technology.

In this way, to $x = \cos \alpha$ corresponds the family of functions $x = \text{cex } \theta = \text{cex } (\theta, s, \varepsilon)$ where $s = e/R$ and $\varepsilon$ are the polar coordinates of the **ex-center $S(s,\varepsilon)$**, which corresponds to the unity/trigonometric circle or **E(e, ε),** which corresponds to a certain circle of radius R, considered as **pole** of a straight line **d,** which rotates around **E** or **S** with the position angle **θ,** generating in this way the ex-centric trigonometric functions, or ex-centric circular supermathematics functions (**EC-SMF**), by intersecting **d** with the unity circle (see.Fig.**1**). Amongst them the **ex-centric cosine** of **θ,** denoted **cex θ** = **x**, where **x** is the projection of the point **W**, which is the intersection of the straight line with the trigonometric circle **C(1,O),** or the Cartesian coordinates of the point **W**. Because a straight line, passing through **S**, interior to the circle (**s ≤ 1 → e < R**), intersects the circle in two points **W₁** and **W₂**, which can be denoted **W$_{1,2}$**, it results that there are **two determinations** of the ex-centric circular supermathematics functions (EC-SMF): a principal one of index **1 cex$_1$ θ,** and a secondary one **cex$_2$ θ**, of index **2**, denoted **cex$_{1,2}$ θ**. **E** and **S** were named **ex-centre** because they were excluded from the center **O(0,0).** This exclusion leads to the apparition of EM and implicitly of SM. By this, the number of mathematical objects grew from one to infinity: to a **unique** function from CM, for example **cos α,** corresponds an **infinity** of functions **cex θ**, due to the possibilities of placing the **ex-center S** and/or **E** in the plane.

**S(e, ε)** can take an infinite number of positions in the plane containing the unity or trigonometric circle. For each position of **S** and **E** we obtain a function **cex θ.** If **S** is a fixed point, then we obtain the ex-centric circular SM functions (EC-SMF), with fixed **ex-center**, or with constant **s** and **ε**. But **S** or **E** can take different positions, in the plane, by various rules or laws, while the straight line which generates the functions by its intersection with the circle, rotates with the angle **θ** around **S** and **E.**

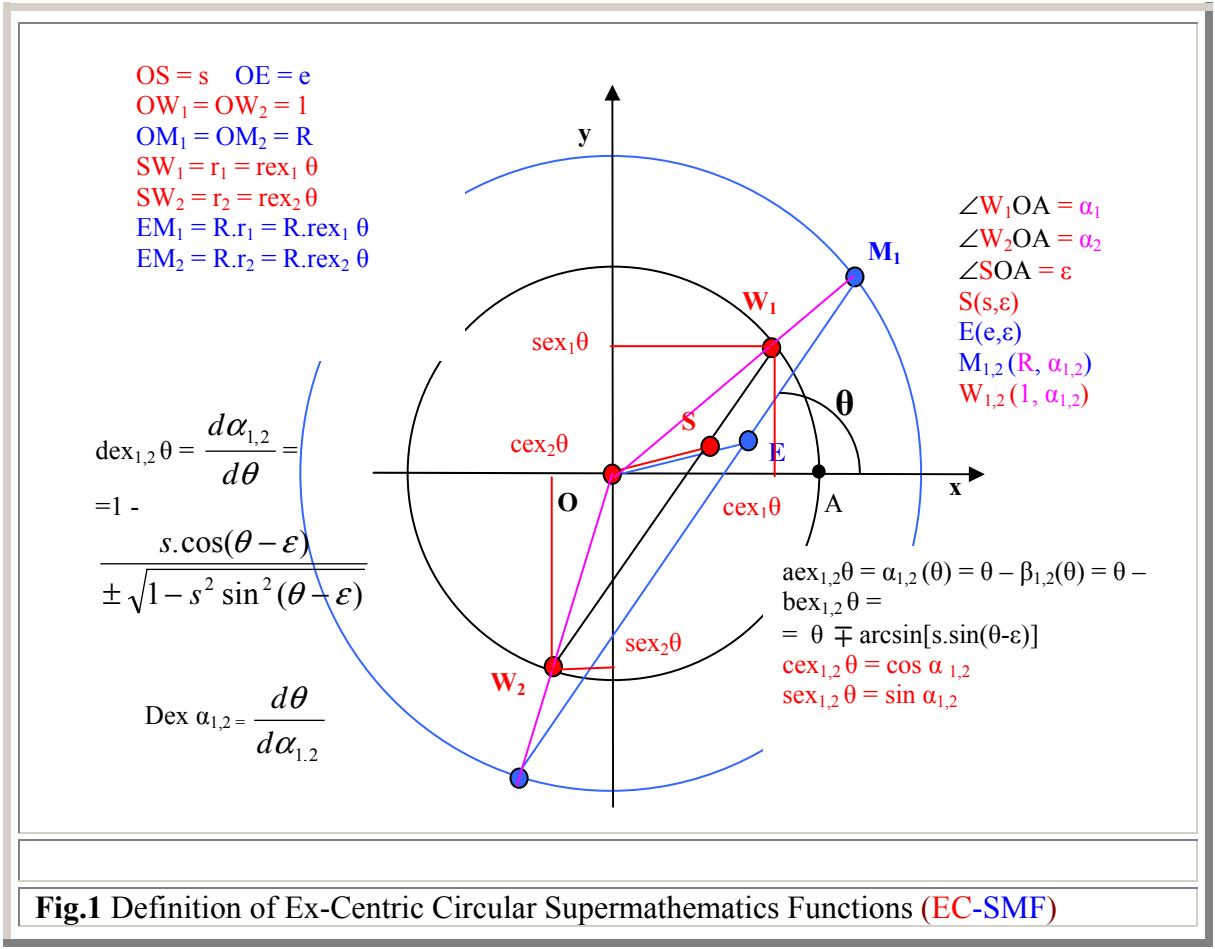

**Fig.1** Definition of Ex-Centric Circular Supermathematics Functions (EC-SMF)

In the last case, we have an EC-SMF of ex-center variable point S/E, which means **s = s (θ)** and/or **ε = ε (θ).** If the variable position of **S/E** is represented also by EC-SMF of the same ex-center **S(s, ε)** or by another ex-center $S_1[s_1 = s_1(θ), ε_1 = ε_1(θ)]$, then we obtain functions of double ex-centricity. By extrapolation, we'll obtain functions of triple, and multiple ex-centricity. Therefore, **EC**-SMF are functions of as many variables as we want or as many as we need**.**

If the distances from **O** to the points $W_{1,2}$ on the circle **C(1,O)** are constant and equal to the radius **R = 1** of the trigonometric circle **C**, distances that will be named **ex-centric radiuses**, the distances from **S** to $W_{1,2}$ denoted by $r_{1,2}$ are variable and are named ex-centric radiuses of the unity circle **C(1,O)** and represent, in the same time, new ex-centric circular supermathematics functions (EC-SMF), which were named **ex-centric radial functions,** denoted $rex_{1,2}\ θ$, if are expressed in function of the **variable** named **ex-centric θ** and **motor**, which is the angle from the ex-center **E**. Or, denoted $Rex1,2\ α,$ if it is expressed in function of the angle **α** or the **centric variable**, the angle at **O(0,0)**. The $W_{1,2}$ are seen under the angles $α_{1,2}$ from **O(0,0)** and under the angles **θ** and **θ + π** from **S(e, ε)** and **E**. The straight line **d** is divided by **S ⊂ d** in the two semi-straight lines, one positive **d $^+$** and the other negative **d $^-$** . For this reason, we can consider $r_1$ = $rex_1$ **θ** a positive oriented segment on **d** (→ **r₁ > 0**) and **r₂** = rex₂ **θ** a negative oriented segment on **d** (→ **r₂ < 0**) in the negative sense of the semi-straight line **d⁻**.

Using simple trigonometric relations, in certain triangles **OEW₁,₂**, or, more precisely, writing the sine theorem (as function of **θ**) and Pitagora's generalized theorem (for the variables **α₁,₂**) in these triangles, it immediately results the **invariant expressions** of the ex-centric radial functions:

$$r_{1,2}(\theta) = rex_{1,2}\,\theta = -s\cos(\theta - \varepsilon) \pm \sqrt{1 - s^2 \sin^2(\theta - \varepsilon)}$$

and

$$r_{1,2}(\alpha_{1,2}) = Rex\,\alpha_{1,2} = \pm\sqrt{1 + s^2 - 2s\cos(\theta - \varepsilon)}.$$

All **EC-SMF** have **invariant** expressions, and because of that they don't need to be tabulated, tabulated being only the centric functions from **CM**, which are used to express them. In all of their expressions, we will always find one of the square roots of the previous expressions, of ex-centric radial functions.

Finding these two determinations is simple: for **+ (plus)** in front of the square roots we always obtain the first determination (**r₁ > 0**) and for the **─ (minus)** sign we obtain the second determination (**r₂ < 0**). The rule remains true for all EC-SMF. By convention, the first determination, of index **1**, can be used or written without index.

Some remarks about these REX ("King") functions:
- The ex-centric radial functions are the expression of the distance between two points, in the plane, in polar coordinates: $S(s,\varepsilon)$ and $W_{1,2}$ ($R = 1$, $\alpha_{1,2}$), on the direction of the straight line **d**, skewed at an angle **θ** in relation to **Ox** axis;
- Therefore, using exclusively these functions, we can express the equations of all known **plane curves,** as well as of other new ones, which surfaced with the introduction of **EM.** An example is represented by **Booth's lemniscates** (see Fig. **2, a, b, c**), expressed, in polar coordinates, by the equation:
  $\rho(\theta) = R(rex_1\theta + rex_2\theta) = -2\,s.R\cos(\theta - \varepsilon)$ for **R=1**, **ε = 0** and **s ∈ [0, 3]**.

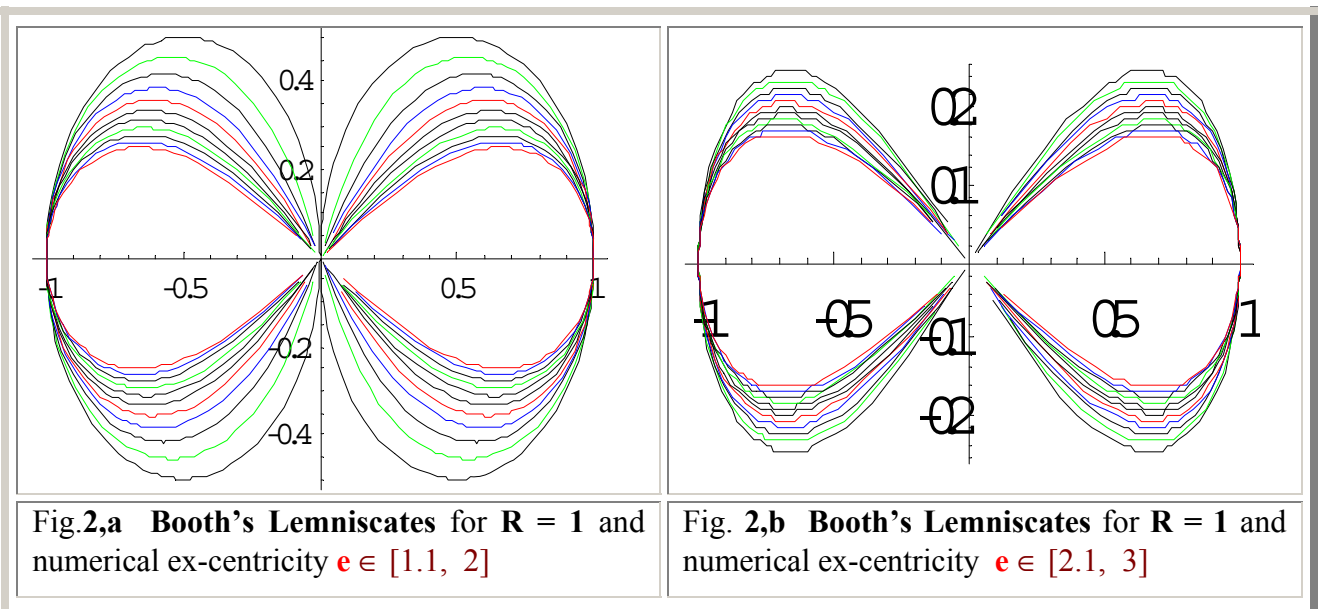

| Fig.**2,a** Booth's Lemniscates for **R = 1** and numerical ex-centricity **e ∈ [1.1, 2]** | Fig. **2,b** Booth's Lemniscates for **R = 1** and numerical ex-centricity **e ∈ [2.1, 3]** |

- Another consequence is the generalization of the definition of a circle:

"**The Circle** is the plane curve whose points M are at the distances **r(θ) = R.rex θ = R.rex [θ, E(e, ε)]** in relation to a **certain** point from the circle's plane **E(e, ε)**".
If **S ≡ O(0,0)**, then **s = 0** and rex θ = 1 = constant, and **r(θ) = R** = constant, we obtain the circle's **classical definition**: the points situated at the same distance **R** from a point, the center of the circle.

**Booth Lemniscate Functions**

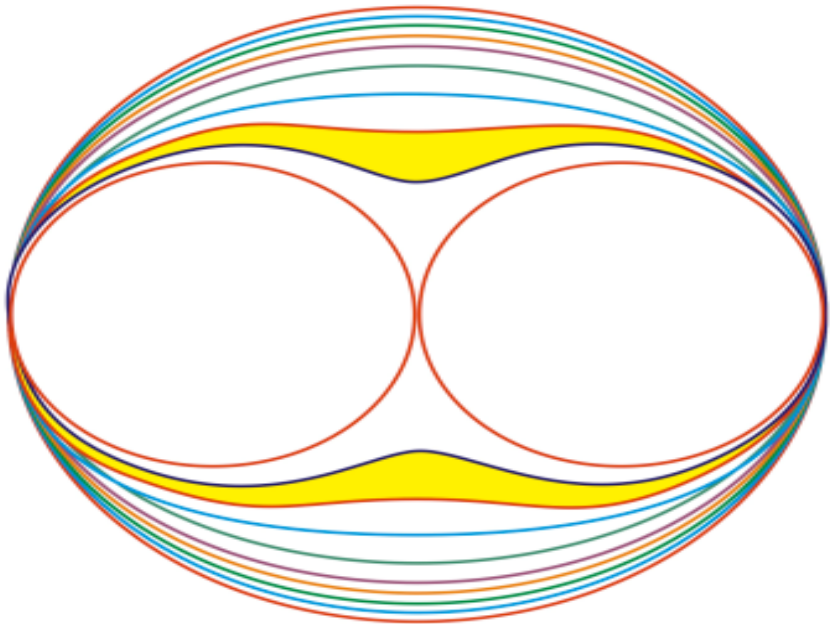

**Polar coordinate equation with**
**supermathematics circle functions rex $_{1,2}$ θ :**
**ρ = R (rex$_1$ θ + rex$_2$ θ )**
**for**
**circle radius R = 1**
**and**
**the numerical ex-centricity s ∈ [ 0, 1 ]**

Fig. **2,c**

- The functions **rex θ** and **Rex α** expresses the transfer functions of zero degree, or of the position of transfer, from the mechanism theory, and it is the ratio between the parameter **R(α$_{1,2}$)**, which positions the conducted element **OM$_{1,2}$** and parameter **R.r$_{1,2}$(θ),** which positions the leader element **EM$_{1,2}$**.
  Between these two parameters, there are the following relations, which can be deduced similarly easy from Fig. 1 that defines EC-SMF.
  Between the position angles of the two elements, leaded and leader, there are the following relations:
  $$α_{1,2} = θ \Upsilon \arcsin[e.\sin(θ − ε)] = θ \Upsilon β_{1,2}(θ) = aex_{1,2} θ$$

  and

  $$θ = α_{1,2} \pm β_{1,2}(α_{1,2}) = α_{1,2} \pm \arcsin\left[ \pm \frac{s.\sin(α_{1,2} - ε)}{\sqrt{1 + s^2 - 2.s.\cos(α_{1,2} - ε)}} \right] = Aex(α_{1,2}).$$

  The functions **aex $_{1,2}$ θ** and **Aex α$_{1,2}$** are EC-SMF, called **ex-centric amplitude**, because of their usage in defining the ex-centric cosine and sine from EC-SMF, in the same manner as the amplitude function or amplitudinus am(k,u) is used for defining the elliptical Jacobi functions:
  sn(k,u) = sn[am(k,u)], cn(k,u) = cos[am(k,u)],
  or:

  **cex$_{1,2}$ θ = cos(aex$_{1,2}$ θ)**,         **Cex α$_{1,2}$ = cos(Aex α$_{1,2}$)**
  and
  **sex $_{1,2}$ θ = sin (aex$_{1,2}$ θ),**        **Sex α$_{1,2}$ = cos (Aex α$_{1,2}$)**

- The radial ex-centric functions can be considered as modules of the position vectors $\vec{r}_{1,2}$ for the **W$_{1,2}$** on the unity circle **C (1,O)**. These vectors are expressed by the following relations:
  $$\vec{r}_{1,2} = rex_{1,2}θ.rad\,θ,$$

  where **rad θ** is the unity vector of variable direction, or the versor/**phasor** of the straight line direction **d$^+$**, whose derivative is the phasor **der θ = d(rad θ)/d θ** and represents normal vectors on the straight lines **OW$_{1,2,}$** directions, tangent to the circle in the **W$_{1,2}$**. They are named the **centric derivative** phasors. In the same time, the modulus **rad θ** function is the corresponding, in **CM,** of the function **rex θ** for **s = 0 →  θ = α** when **rex θ = 1** and **der α$_{1,2}$** are the tangent versors to the unity circle in **W$_{1,2}$**.

- The derivative of the $\vec{r}_{1,2}$ vectors are the velocity vectors:
  $$\vec{v}_{1,2} = \frac{d\vec{r}_{1,2}}{dθ} = dex_{1,2}θ.der\,α_{1,2}$$
  of the **W$_{1,2}$ ⊂ C** points in their rotating motion on the circle, with velocities of variable modulus **v$_{1,2}$ = dex$_{1,2}$ θ**, when the generating straight line **d** rotates around the ex-center **S** with a constant angular speed and equal to the unity, namely **Ω = 1.** The velocity vectors have the expressions presented above, where **der** $\alpha_{1,2}$ are the phasors of centric radiuses $R_{1,2}$ of module **1** and of $\alpha_{1,2}$ directions. The expressions of the functions **EC-SM** **dex**$_{1,2}$ $\theta$, **ex-centric derivative** of $\theta$, are, in the same time, also the $\alpha_{1,2}(\theta)$ angles derivatives, as function of the motor or independent variable $\theta$, namely

$$\text{dex}_{1,2}\,\theta = d\alpha_{1,2}(\theta)/d\theta = 1 - \frac{s.\cos(\theta-\varepsilon)}{\pm\sqrt{1-s^2.\sin^2(\theta-\varepsilon)}}$$

as function of $\theta$, and

$$\textbf{Dex}\,\alpha_{1,2} = d(\theta)/d\alpha_{1,2} = \frac{1-s.\cos(\alpha_{1,2}-\varepsilon)}{1+s^2 - 2.s.\cos(\alpha_{1,2}-\varepsilon)} = \frac{1-s.\cos(\alpha_{1,2}-\varepsilon)}{\text{Re}\,x^2\alpha_{1,2}},$$

as functions of $\alpha_{1,2}$.

It has been demonstrated that the **ex-centric derivative** functions **EC-SM** express the transfer functions of the first order, or of the angular velocity, from the Mechanisms Theory, for **all (!)** known plane mechanisms.

- The radial ex-centric function **rex** $\theta$ expresses exactly the movement of push-pull mechanism **S = R. rex** $\theta$, whose motor connecting rod has the length **r**, equal with **e** the real ex-centricity, and the length of the crank is equal to **R**, the radius of the circle, a very well-known mechanism, because it is a component of all automobiles, except those with Wankel engine.

The applications of radial ex-centric functions could continue, but we will concentrate now on the more general applications of **EC-SMF**.

Concretely, to the unique forms as those of the circle, square, parabola, ellipse, hyperbola, different spirals, etc. from **CM**, which are now grouped under the name of **centrics**, correspond an infinity of **ex-centrics** of the same type**:** circular, square (quadrilobe), parabolic, elliptic, hyperbolic, various spirals **ex-centrics**, etc. Any **ex-centric function**, with null ex-centricity (e = 0), degenerates into a **centric function**, which represents, at the same time its **generating** curve. Therefore, the CM itself belongs to **EM**, for the unique case (s = e = 0), which is one case from an infinity of possible cases, in which a point named eccenter **E(e, ε)** can be placed in plane. In this case, **E** is overleaping on one or two points named **center**: the **origin O(0,0)** of a frame, considered the origin **O(0,0)** of the referential system, and/or the **center C(0,0)** of the unity circle for circular functions, respectively, the symmetry center of the two arms of the equilateral hyperbola, for hyperbolic functions.

It was enough that a point **E** be eliminated from the center (**O** and/or **C**) to generate from the old CM a new world of **EM**. The reunion of these two worlds gave birth to the **SM** world.

This discovery occurred in the city of the Romanian Revolution from 1989, **Timişoara,** which is the same city where on November 3[rd], 1823 **Janos Bolyai** wrote: "From **nothing** I've created a new world". With these words, he announced the discovery of the fundamental formula of the first **non-Euclidean geometry**.

He – from nothing, I – in a joint effort, proliferated the periodical functions which are so helpful to engineers to describe some periodical phenomena. In this way, I have enriched the mathematics with new objects.

When **Euler** defined the trigonometric functions, as direct circular functions, if he wouldn't have chosen **three superposed points**: the origin **O**, the center of the circle **C** and **S** as a pole of a semi straight line, with which he intersected the trigonometric/unity circle, the **EC-SMF** would have been discovered much earlier, eventually under another name.

Depending on the way of the "**split**" (we isolate one point at the time from the superposed ones, or all of them at once), we obtain the following types of **SMF**:

**O ≡ C ≡ S** → **Centric functions belonging to CM;**

and those which belong to **EM** are:

**O ≡ C ≠ S** → **Ex-centric Circular Supermathematics Functions (EC-SMF);**
**O ≠ C ≡ S** → **Elevated Circular Supermathematics Functions (ELC-SMF);**
**O ≠ C ≠ S** → **Exotic Circular Supermathematics Functions (EXC-SMF).**

These **new mathematics complements**, joined under the **temporary** name of **SM**, are extremely useful tools or instruments, long awaited for. The proof is in the large number and the diversity of periodical functions introduced in mathematics, and, sometimes, the complex way of reaching them, by trying the substitution of the circle with other curves, most of them closed.

To obtain new special, periodical functions, it has been attempted the replacement of the trigonometric circle with the **square** or the **diamond**. This was the proceeding of Prof. Dr. Math. **Valeriu Alaci,** the former head of the Mathematics Department of Mechanics College from Timişoara, who discovered the **square** and **diamond** trigonometric functions. Hereafter, the mathematics teacher **Eugen Visa** introduced the **pseudo-hyperbolic** functions, and the mathematics teacher **M. O**. **Enculescu** defined the **polygonal** functions, replacing the circle with an n-sides polygon; for n = 4 he obtained the square **Alaci** trigonometric functions. Recently, the mathematician, Prof. **Malvina Baica,** (of Romanian origin) from the University of Wisconsin together with Prof**. Mircea Cárdu,** have completed the gap between the **Euler** circular functions and **Alaci** square functions, with the so-called **Periodic Transtrigonometric functios.**

The Russian mathematician **Marcusevici** describes, in his work "**Remarkable sine functions**" the **generalized trigonometric functions** and the trigonometric functions **lemniscates**.

Even since 1877, the German mathematician Dr. **Biehringer**, substituting the right triangle with an oblique triangle, has defined the **inclined** trigonometric functions. The British scientist of Romanian origin Engineer **George** (Gogu) **Constantinescu** replaced the **circle** with the **evolvent** and defined the **Romanian** trigonometric functions: **Romanian cosine** and **Romanian sine**, expressed by Cor α and Sir α functions, which helped him to resolve some non-linear differential equations of the Sonicity Theory, which he created. And how little known are all these functions even in Romania!

Also the elliptical functions are defined on an ellipse. A rotated one, with its main axis along Oy axis.

How simple the complicated things can become, and as a matter of fact they are! This **paradox**(ism) suggests that by a simple displacement/expulsion of a **point** from a **center** and by the apparition of the notion of the **ex-center**, a new world appeared, the world of EM and, at the same time, a new Universe, the **SM** Universe.

Notions like "**Supermathematics Functions**" and "**Circular Ex-centric Functions**" appeared on most search engines like Google, Yahoo, AltaVista etc., from the beginning of the Internet. The new notions, like **quadrilobe** "quadrilobas", how the **ex-centrics** are named, and which continuously fill the space between a **square** circumscribed to a circle and the **circle** itself were included in the Mathematics Dictionary. The intersection of the **quadriloba** with the straight line **d** generates the new functions called **cosine quadrilobe**-ic and **sine quadrilobe**-ic**.**

The benefits of **SM** in science and technology are too numerous to list them all here. But we are pleased to remark that **SM** removes the boundaries between **linear** and **non-linear**; the linear belongs to **CM**, and the non-linear is the appanage of **EM**, as between **ideal** and **real**, or as between perfection and **imperfection**.

It is known that the **Topology** does not differentiate between a pretzel and a cup of tea. Well, **SM** does not differentiate between a **circle** (e = 0) and a **perfect square** (s = ± 1), between a **circle** and a **perfect triangle**, between an **ellipse** and a **perfect rectangle**, between a **sphere** and a perfect **cube,** etc. With the same parametric equations we can obtain, besides the **ideal** forms of **CM** (circle, ellipse, sphere etc.), also the **real** ones (square, oblong, cube, etc.). For **s** ∈ [-1,1], in the case of ex-centric functions of variable **θ**, as in the case of centric functions of variable **α,** for **s**∈[-∞,+∞], it can be obtained an infinity of intermediate forms, for example, square, oblong or cube with rounded corners and slightly curved sides or, respectively, faces. All of these facilitate the utilization of the new SM functions for drawing and representing of some technical parts, with rounded or splayed edges, in the **CAD/ CAM-SM** programs, which don't use the computer as drawing boards any more, but create the technical object instantly, by using the parametric equations, that speed up the processing, because only the equations are memorized, not the vast number of pixels which define the technical piece.

The numerous functions presented here, are introduced in mathematics for the first time, therefore, for a better understanding, the author considered that it was necessary to have a short presentation of their equations, such that the readers, who wish to use them in their application's development, be able to do it.

**SM** is not a finished work; it's merely an **introduction** in this vast domain, a first step, the author's small step, and a giant leap for mathematics.

The **elevated circular SM** functions (**ELC-SMF**), named this way because by the modification of the numerical ex-centricity **s** the points of the curves of elevated sine functions sel **θ** as of the elevated circular function elevated cosine **cel θ** is elevating – in other words it rises on the vertical, getting out from the space {-1, +1} of the other sine and cosine functions, centric or ex-centric. The functions' **cex θ** and **sex θ** plots are shown in Fig. **3,** where it can be seen that the points of these graphs get modified on the horizontal direction, but all remaining in the space
[-1,+1], named the existence domain of these functions.

The functions' cel **θ** and sel **θ** plots can be simply represented by the products:

$$\text{cel}_{1,2} \ \theta = \text{rex}_{1,2} \ \theta \cdot \cos \theta \quad \text{and} \quad \text{Cel } \alpha_{1,2} = \text{Rex } \alpha_{1,2} \cdot \cos \theta$$
$$\text{sel}_{1,2} \ \theta = \text{rex}_{1,2} \ \theta \cdot \sin \theta \quad \text{and} \quad \text{Sel } \alpha_{1,2} = \text{Rex } \alpha_{1,2} \cdot \sin \theta$$

and are shown Fig. **4.**

The **exotic circular functions** are the most general **SM,** and are defined on the unity circle which is not centered in the origin of the xOy axis system, neither in the eccenter **S,** but in a certain point **C (c,γ)** from the plane of the unity circle, of polar coordinates (**c, γ**) in the xOy coordinate system.

Many of the drawings from this album are done with **EC-SMF** of ex-center variable and with arcs that are multiples of **θ (n.θ).** The used relations for each particular case are explicitly shown, in most cases using the **centric** mathematical functions, with which, as we saw, we could express all **SM** functions, especially when the image programs cannot use SMF. This doesn't mean that, in the future, the new math complements will not be implemented in computers, to facilitate their vast utilization.

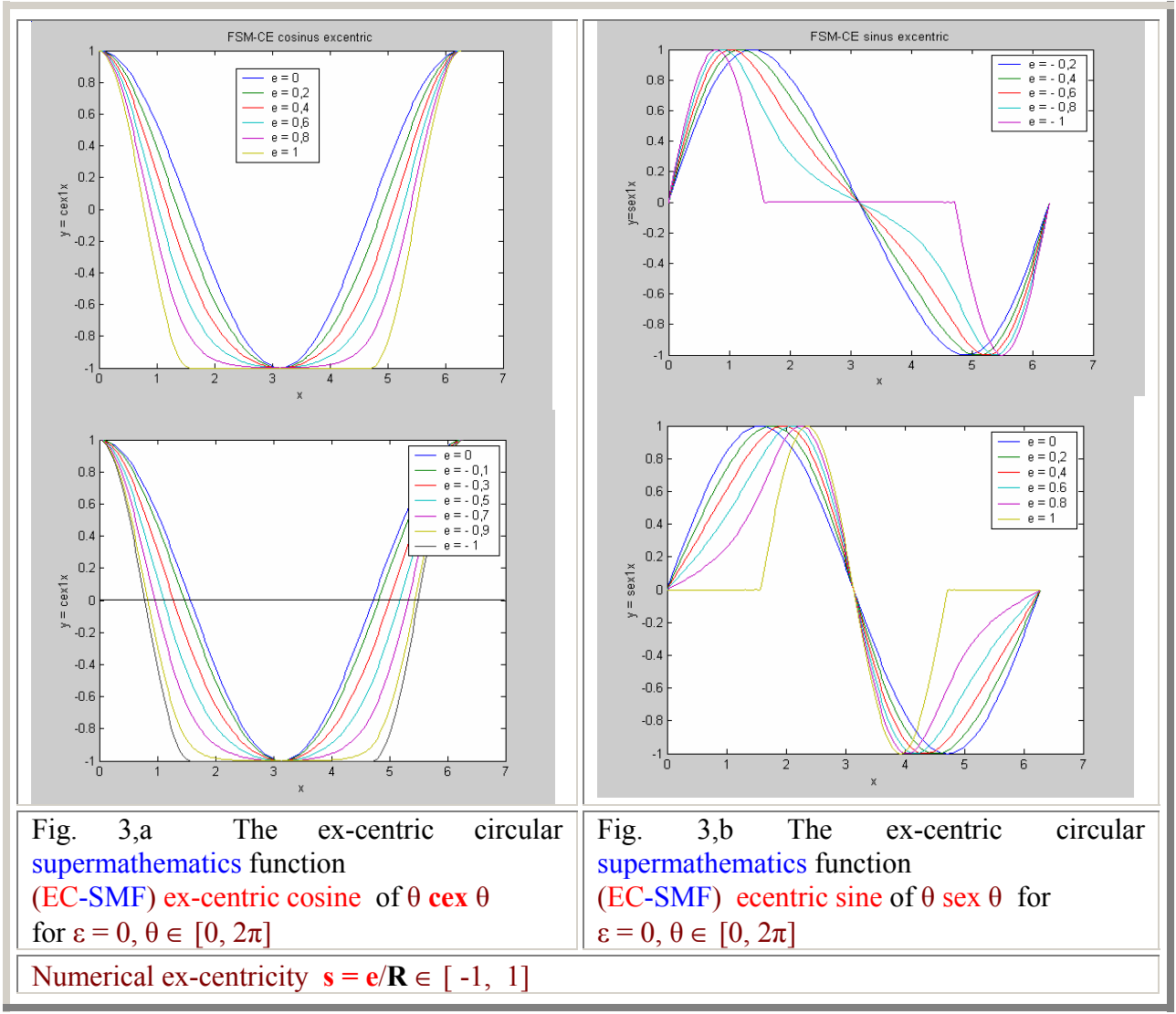

| Fig. 3,a The ex-centric circular supermathematics function (EC-SMF) ex-centric cosine of θ **cex** θ for ε = 0, θ ∈ [0, 2π] | Fig. 3,b The ex-centric circular supermathematics function (EC-SMF) ecentric sine of θ **sex** θ for ε = 0, θ ∈ [0, 2π] |
|---|---|
| Numerical ex-centricity **s = e/R** ∈ [ -1, 1] ||

The computer specialists working in programming the computer assisted design software **CAD/CAM/CAE,** are on their way to develop these new programs fundamentally different, because the technical objects are created with parametric **circular** or **hyperbolic** SMFs, as it has been exemplified already with some achievements such as airplanes, buildings, etc. in http://www.eng.upt.ro/~mselariu and how a washer can be represented as a toroid ex-centricity (or as an "ex-centric torus"), square or oblong in an axial section, and, respectively, a square plate with a central square hole can be a "square torus of square section". And all of these, because **SM** doesn't make distinction between a circle and a square or between an ellipse and a rectangle, as we mentioned before.

But the most important achievements in science can be obtained by solving some non-linear problems, because **SM** reunites these two domains, so different in the past, in a single entity. Among these differences we mention that the non-linear domain asks for ingenious approaches for each problem. For example, in the domain of vibrations, static elastic characteristics (SEC) soft non-linear (regressive) or hard non-linear (progressive) can be obtained simply by writing y = m. x, where **m** is not anymore m = tan α as in the linear case  (**s** = 0 ), but  **m** = tex$_{1,2}$ **θ** and depending on the numerical ex-centricity **s** sign, positive or negative, or for **S** placed on the negative x axis (ε = π) or on the positive x axis (ε = 0), we obtain the two nonlinear elastic characteristics, and obviously for s=0  we'll obtain the linear SEC.

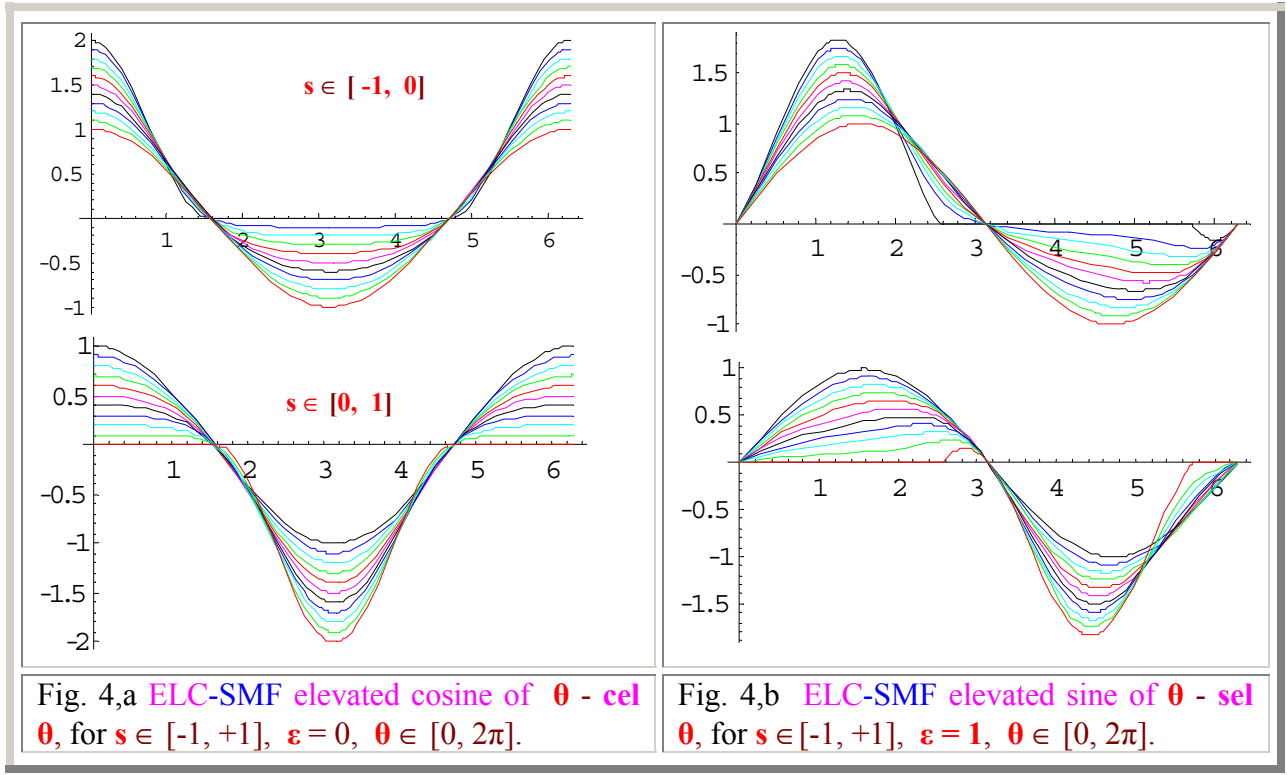

| Fig. 4,a ELC-SMF elevated cosine of  θ - c**el** θ, for **s** ∈ [-1, +1], ε = 0, θ ∈ [0, 2π]. | Fig. 4,b  ELC-SMF elevated sine of θ - s**el** θ, for **s** ∈[-1, +1],  ε = 1,  θ ∈ [0, 2π]. |

Due to the fact that the functions **cex θ** and **sex θ**, as well **Cex** α and **Sex α** and their combinations, are solutions of some differential equations of second degree with variable coefficients, it has been stated that the linear systems (Tchebychev) are obtained also for **s = ± 1,** and not only for **s = 0.** In these equations, the mass ( the point **M**) rotates on the circle with a double angular speed  **ω = 2.Ω**  (reported to the linear system where  **s = 0** and **ω = Ω = constant**) in a half of a period, and in the other half of period stops in the point **A(R**,0)  for  **e = sR = R**  or  ε = 0 and in **A'(−R**, 0) for  **e = − s.R = −1**, or  ε = π. Therefore, the oscillation period **T** of the **three linear systems** is the same and equal with **T = Ω / 2π**.  The nonlinear S**E**C systems are obtained for the others values, intermediates, of **s** and **e**. The projection, on any direction, of the rotating motion of **M** on the circle with radius **R**, equal to the oscillation amplitude, of a variable angular speed **ω = Ω.dex θ**  ( after **dex θ** function) is an **non-linear**  oscillating motion**.**

The discovery of  "**king**" function **rex θ,** with its properties, facilitated the apparition of a **hybrid** method (analytic-numerical), by which a simple **relation** was obtained, with only two terms, to **calculate** the first degree elliptic complete integral  **K(k)**, with an unbelievable precision, with a minimum of **15 accurate decimals,** after only 5 steps. Continuing with the next steps, can lead us to a new relation to compute **K (k)**, with a considerable higher precision and with possibilities to expand the method to other elliptic integrals, and not only to those. After 6 steps, the relation of **E (k)** has the same precision of computation.

The discovery of **SMF** facilitated the apparition of a new integration method, named **integration through the differential dividing**.

We will stop here, letting to the readers the pleasure to delight themselves by viewing the drawings from this album.

**References** (in SuperMathematics):

| 13 | Selariu Mircea | RIGIDITATEA DINAMICA EXPRIMATA CU FUNCTII SUPERMATEMATICE | Com.VII Conf. Internat. de Ing. Manag. si Tehn., TEHNO'95 Timisoara, 1995 Vol.7: Mecatronica, Dispoz. si Rob.Ind.,pag. 185...194 |
|---|---|---|---|
| 14 | Selariu Mircea | DETERMINAREA ORICAT DE EXACTA A RELATIRI DE CALCUL A INTEGRALEI ELIPTICE COMPLETE D E SPETA INTAIA K(k) | Bul. VIII-a Conf. de Vibr. Mec., Timisoara,1996, Vol III, pag.15 ... 24. |
| 15 | Selariu Mircea | FUNCTII SUPERMATEMATICE CIRCULARE EXCENTRICE DE VARIABILA CENTRICA | A VIII_a Conf. Internat. de Ing. Manag. si Tehn. TEHNO'98, Timisoara, 1998, pag. 531...548 |
| 16 | Selariu Mircea | FUNCTII DE TRANZITIE INFORMATIONALA | A VIII_a Conf. Internat. de Ing. Manag. si Tehn. TEHNO'98, Timisoara, 1998, pag.549..556 |
| 17 | Selariu Mircea | FUNCTII SUPERMATEMATICE EXCENTRICE DE VARIABILA CENTRICA CA SOLUTII ALE UNOR SISTEME OSCILANTE NELINIARE | A VIII_a Conf. Internat. de Ing. Manag. si Tehn. TEHNO'98, Timisoara, 1998, pag. 557..572 |
| 18 | Selariu Mircea | TRANSFORMAREA RIGUROASA IN CERC A DIAGRAMEI POLARE A COMPLIANTEI | Bul. X Conf. VCM ,Bul St. Si Tehn. Al Univ. Poli. Timisoara, Seria Mec. Tom. 47 (61) mai 2002, Vol II pag. 247…260 |
| 19 | Selariu Mircea | INTRODUCEREA STRAMBEI IN MATEMATICA | Luc.Simp. Nat. Al Univ. Gh. Anghel Drobeta Tr. Severin, mai 2003, pag. 171…178 |
| 20 | Petrisor Emilia | ON THE DYNAMICS OF THE DEFORMED STANDARD MAP | Workshop Dynamicas Days'94, Budapest, si Analele Univ.din Timisoara, Vol.XXXIII, Fasc.1-1995, Seria Mat.-Inf.,pag. 91…105 |
| 21 | Petrisor Emilia | SISTEME DINAMICE HAOTICE | Seria Monografii matematice, Tipografia Univ. de Vest din Timisoara, 1992 |
| 22 | Petrisor Emilia | RECONNECTION SCENARIOS AND THE THERESHOLD OF RECONNECTION IN THE DYNAMICS OF NONTWIST MAPS | Chaos, Solitons and Fractals, 14 ( 2002) 117…127 |
| 23 | Cioara Romeo | FORME CLASICE PENTRU FUNCTII CIRCULARE EXCENTRICE | Proceedings of the Scientific Communications Meetings of "Aurel Vlaicu" University, Third Edition, Arad, 1996, pg.61 ...65 |
| 24 | Preda Horea | REPREZENTAREA ASISTATA A TRAIECTORILOR IN PLANUL FAZELOR A | Com. VI-a Conf.Nat.Vibr. in C.M. Timisoara, 1993, pag. |

| | | VIBRATIILOR NELINIARE | |
|---|---|---|---|
| 25 | Selariu Mircea Ajiduah Crist. Bozantan Emil (USA) Filipescu Avr. | INTEGRALELE UNOR FUNCTII SUPERMATEMATICE | Com. VII Conf.Intern.de Ing.Manag.si Tehn. TEHNO'95 Timisoara. 1995,Vol.IX: Matem.Aplic. pag.73...82 |
| 26 | Selariu Mircea | CALITATEA CONTROLULUI CALITATII | Buletin AGIR anul II nr.2 (4) -1997 |
| 27 | Selariu Mircea Fritz Georg (G) Meszaros A. (G) | ANALIZA CALITATII MISCARILOR PROGRAMATE cu FUNCTII SUPERMATEMATICE | IDEM, Vol.7: Mecatronica, Dispozitive si Rob.Ind., pag. 163...184 |
| 28 | Selariu Mircea Szekely Barna ( Ungaria ) | ALTALANOS SIKMECHANIZMUSOK FORDULATSZAMAINAK ATVITELI FUGGVENYEI MAGASFOKU MATEMATIKAVAL | Bul.St al Lucr. Prem.,Universitatea din Budapesta, nov. 1992 |
| 29 | Selariu Mircea Popovici Maria | A FELSOFOKU MATEMATIKA ALKALMAZASAI | Bul.St al Lucr. Prem., Universitatea din Budapesta, nov. 1994 |
| 30 | Konig Mariana Selariu Mircea | PROGRAMAREA MISCARII DE CONTURARE A ROBOTILOR INDUSTRIALI cu AJUTORUL FUNCTIILOR TRIGONOMETRICE CIRCULARE EXCENTRICE | MEROTEHNICA, Al V-lea Simp. Nat.de Rob.Ind.cu Part .Internat. Bucuresti, 1985 pag.419...425 |
| 31 | Konig Mariana Selariu Mircea | PROGRAMAREA MISCARII de CONTURARE ale R I cu AJUTORUL FUNCTIILOR TRIGONOMETRICE CIRCULARE EXCENTRICE, | Merotehnica, V-lea Simp. Nat.de RI cu participare internationala, Buc.,1985, pag. 419 ... 425. |
| 32 | Konig Mariana Selariu Mircea | THE STUDY OF THE UNIVERSAL PLUNGER IN CONSOLE USING THE ECCENTRIC CIRCULAR FUNCTIONS | Com. V-a Conf. PUPR, Timisoara, 1986, pag.37...42 |
| 33 | Staicu Florentiu Selariu Mircea | CICLOIDELE EXPRIMATE CU AJUTORUL FUNCTIEI SUPERMATEMATICE REX | Com. VII Conf. Internationala de Ing.Manag. si Tehn ,Timisoara "TEHNO'95"pag.195-204 |
| 34 | Gheorghiu Em. Octav Selariu Mircea Bozantan | FUNCTII CIRCULARE EXCENTRICE DE SUMA si DIFERENTA DE ARCE | Ses.de com.st.stud.,Sectia Matematica,Timisoara, Premiul II pe 1983 |